\documentclass[a4paper,draft]{article}
\usepackage{amsmath}
\newtheorem{theorem}{Theorem}[section]
\newtheorem{lemma}{Lemma}[section]
\newtheorem{corollary}{Corollary}[section]

\newtheorem{proposition}{Proposition}[section]
\newenvironment{remark}{\vspace{1ex}{\bf Remark.}\rm}{\vspace{1ex}}
\newenvironment{proof}{{\bf Proof.}}{\par\hspace{25em}\rule{1ex}{1ex}\vspace{1ex}}
\newcommand{\ds}{\displaystyle}

\title{Totally geodesic property of the Hopf vector
field.\footnote{Acta   Math. Hungar. 101, 1-2 (2003), 73-92. }}
\author{Yampolsky A.}
\date{}
\begin{document}
\maketitle
\begin{abstract}
    We prove that the Hopf vector field is a unique one among geodesic unit vector fields
    on spheres such that the submanifold generated by the field is totally geodesic in the
    unit tangent bundle with Sasaki metric. As application, we give a new proof of stability
    (instability) of the Hopf vector field with respect to volume variation using
    standard approach from the theory of submanifolds and find exact boundaries for
    the sectional curvature of the Hopf vector field.
    \footnote{
    \mbox{\qquad}{\it AMS subject class:} Primary 54C40,14E20; Secondary 46E25, 20C20}
\end{abstract}

\section*{Introduction}
    Let $(M,g)$ be $(n+1)$-dimensional Riemannian manifold  with metric $g$
    and $\xi$ a fixed unit vector field on $M$.
    Consider $\xi$ as a local mapping $\xi : M \to T_1M $. Then the image
    $\xi(M) $  is a submanifold in the unit tangent sphere bundle $T_1M$.
    The Sasaki metric on the tangent bundle $TM$ induces the Riemannian metric
    on $T_1M$ and on $\xi(M)$ as well.
    So, one may use the geometrical properties of this submanifold to determine
    geometrical characteristics of a unit vector field. Namely, a unit vector field $\xi$
    is said to be {\it  minimal} if $ \xi(M)$ is  of minimal volume with
    respect to the induced  metric  \cite{G-Z, GM}; $\xi$ is {\it totally
    geodesic} if $\xi(M)$ is totally geodesic submanifold in $T_1M$. A complete description of totally geodesic vector fields on 2-dimensional manifolds of constant curvature
    has been given in \cite{Ym3}. It was proved that only unit sphere (among non-flat space forms) admits totally geodesic unit vector fields (which are \emph{non-geodesic} ones). In contrast to this case, on a 3-dimensional sphere the Hopf (and therefore, geodesic) unit vector field is a unique one with globally minimal volume \cite{G-Z}. Since totally geodesic submanifold is always minimal, it was natural to suspict that the Hopf vector field is totally geodesic in fact. In addition, for the spheres of greater dimensions the Hopf vector fields stays minimal but become unstable \cite{J,G-Z}. Nevertheless, something should  distinguish the Hopf vector fields. And this \emph{something}, as we prove here, is their totally geodesic property. Namely, ( Theorem \ref{Hopf} )
\vspace{1ex}

\textit{    Let $\xi$ be a unit geodesic vector field on a unit sphere $S^{n+1}$.
    Then $\xi$ is totally geodesic vector field if and only if $n=2m$
    and $\xi$ is a Hopf vector field.\footnote{The result is true with one additional condition (see below)}     }
\vspace{1ex}

Using this result we give another proof of stability or instability of the Hopf vector
fields with respect to volume functional (cf. \cite{GM1})using a standard formula for
the second volume variation from the theory of the submanifolds ( Theorems \ref{Hopf3}
and \ref{Hopfn} ). \vspace{1ex}

\textsl{The Hopf vector field on $S^{2m+1}$ is \emph{\textbf{stable}} for $m=1$ and
\emph{\textbf{unstable}} for $m>1$.} \vspace{1ex}

We also find an exact boundaries for the sectional cutvature of this field (~Theorem
\ref{Cor}~). \vspace{1ex}

\textit{Sectional curvature of $\xi(S^{2m+1})$ varies between $\frac14$ and $\frac54$.
Minimal bound the curvature achieves on $\xi$-tangential lift (\ref{lifts}) of
$\xi$-section and maximal on $\xi$-tangential lift of $\varphi$-section in terms of
contact metric geometry.} \vspace{1ex}

\section{Some preliminaries}

    Let $(M,g)$ be $(n+1)$-dimensional Riemannian manifold  with metric $g$. Denote by
    $\big<\cdot,\cdot\big>$ a scalar product with respect to $g$ and by $\nabla$ the
    Levi-Civita connection on $M$.

    The {\it Sasaki metric} on $TM$ is defined by the following scalar product:
    if $\tilde X,\tilde Y \in TTM$,  then
\begin{equation}
\label{lab5}
        \big<\big< \tilde X,\tilde Y \big>\big>=
        \big<\pi_* \tilde X, \pi_* \tilde Y\big>+\big<K \tilde X,K \tilde Y\big>
\end{equation}
    where $\pi_*:TTM \to TM $ is a differential of projection $\pi:TM \to M $ and
   $K: TTM \to TM$  is the {\it connection map}.

     Let $\xi$ be a unit vector field on $M$. A vector field $\tilde X \in TTM$ is
     tangent to $\xi(M)$ if and only if \cite{Ym1}
     $$
     \tilde X=(\pi_*\tilde X)^h + (\nabla_{\pi_*\tilde X}\xi)^v,
     $$
    where $(\cdot)^h$ and $(\cdot)^v$ mean {\it horizontal} and {\it vertical} lifts
    of fields into the tangent bundle.

    It is well known that $\xi^v$ is a unit normal vector field on $T_1M$.
    Let $X$ be tangent to $M$, then
    $$
    X^{t}=X^v-\big<X,\xi\big>\xi^v,
    $$
    is always tangent to $T_1M$ and is called the {\it tangential lift } of $X$ \cite{Bx-V3}.
    It is easy to  see that
    $$
    \big<\big<X^{t},Y^{t}\big>\big>=\big<X,Y\big>-\big<X,\xi\big>\big<Y,\xi\big>.
    $$
    Introduce, now, a notion of {\it $\xi$-tangential} and {\it $\xi$-normal} lifts with respect to
    given field $\xi.$ We proceed in the following way.

    Introduce a {\it shape operator} $A_\xi$ for the field $\xi$ by
    $$
    A_\xi X=-\nabla_X\xi,
    $$
    where $X$ is arbitrary vector field on $M$.

    Define a {\it conjugate shape operator} $A_\xi^*$ by
    \begin{equation}\label{conj}
    \big<A_\xi^*X,Y\big>=\big<X,A_\xi Y\big>.
    \end{equation}

    Let $X,Y$ be the vector fields on $M$. Define $\xi$-tangential lift $X_\xi^\tau$ and
    $\xi$-normal lift $Y_\xi^\nu$ of $X$ and $Y$ respectively by
    \begin{equation}\label{lifts}
    \ds
    X_\xi^{\tau}=X^h-(A_\xi X)^{t}, \quad Y_\xi^{\nu}=(A^*_\xi Y)^h+Y^{t}.
    \end{equation}

    Then, $X_\xi^\tau$ is evidently tangent to $\xi(M)$. Moreover,
    $$
    \begin{array}{l}
    \big<\big<X^{\tau}_\xi,Y^{\nu}_\xi\big>\big>=\big<\big<X^h,(A^*_\xi Y)^h\big>\big>-
    \big<\big<(A_\xi X)^{t},Y^{t}\big>\big>= \\[2ex]
     \big<X,A^*_\xi Y\big>-\big<A_\xi X,Y\big>+\big<A_\xi X,\xi\big>\big<Y,\xi\big>=
    \big<A_\xi X,Y\big>-\big<A_\xi X,Y\big>\equiv0
    \end{array}
    $$
    and therefore the linear space of all $\xi$-normal lifts coincides with the normal
    space of $\xi(M)$ at each point.

    To construct a natural tangent and normal orthonormal frames for $\xi(M)$, one can
    use a singular decomposition of the shape operator $A_\xi$, based on the
    following linear algebra result ( \cite{HG}, Theorem 7.3.5 and Exercise 7.3.5)
    which we present here in a slightly modified form.
    \begin{theorem}\label{Sing}
    A matrix  $A\in M_{m,n}$ of rank $k$ may be represented in the form
    $$
    A=F\Sigma E^*,
    $$
    where $F\in M_m$ and $E\in M_n$ are unitary matrices.

    The matrix
    $\Sigma=[\sigma_{ij}]\in M_{m,n}$ is such that
    $\sigma_{ij}=0,\,i\ne j,\,\sigma_{11}\geq\sigma_{22}\geq\dots\geq\sigma_{kk}>
    \sigma_{k+1,k+1}=\dots=\sigma_{qq}=0, \, q=\min\{m,n\}.$ The values
    $\{\sigma_{ii}\}\equiv\{\lambda_i\}$ are non-negative square roots of eigenvalues
    of the matrix $AA^*$ and hence are uniquely defined. The columns of  the matrix
    $F$ are the eigenvectors of the matrix $AA^*$ and the columns of the matrix $E$
    are the eigenvalues of the matrix $A^*A$. Moreover, $A^*f_i=\lambda_i e_i$ and
    $A e_i=\lambda_i f_i$ for $i=1,\dots,k$. If the matrix $A$ is real, then
    $F,\Sigma$ and $E$ can be real.
    \end{theorem}

    The columns of the matrices $F$ and $E$ are called respectively {\it left} and
    {\it right singular vectors } of matrix $A$. The values $\lambda_i$ are called {\it
    singular values } of the matrix $A$.

    Set $A=A_\xi$ and apply Theorem \ref{Sing}.
    Since $A_\xi^*\xi=0$ for any unit vector field $\xi$ ,
    there exist an orthonormal local frames $e_0, e_1, \dots , e_n $ and
    $f_0=\xi,\, f_1, \dots , f_n $ on $M$ such that
    $$
        A_\xi \, e_0=0 , \ A_\xi\, e_\alpha =\lambda_\alpha f_\alpha,
        \quad A_\xi^*\, f_0=0,\ A_\xi^*\,f_\alpha=\lambda_\alpha e_\alpha ,
        \ \alpha=1, \dots , n ,
    $$
    where $\lambda_1\geq \lambda_{2}\geq \dots  \lambda_n\geq0$ are the real-valued
    functions.

    It is natural to call the functions $\lambda_i$ ($ i=1,\dots, n$) the {\it singular principal
    curvatures } of the field $\xi$ with respect to chosen singular frame.
    Remark, that if necessary one may use the signed singular values fixing the directions of the
    vectors of the singular frame. Setting $\lambda_0=0$, we may rewrite the relations
    on singular frames in a unified form
    \begin{equation}\label{sing}
    \begin{array}{l}
     \ds   A_\xi\, e_i =\lambda_i f_i,\quad
        A_\xi^*\,f_i=\lambda_i e_i ,
        \qquad i=0,1, \dots , n ,\\[2ex]
     \ds   \lambda_0=0, \ \lambda_1,\dots,\lambda_n\geq0.
    \end{array}
    \end{equation}

    The following lemma is easy to prove using (\ref{conj}) and (\ref{sing}).
\begin{lemma}\cite{Ym1}
    At each point of $\xi(M)\subset TM$ the orthonormal frames
    \begin{equation}\label{tang}
    \begin{array}{ll}
     \ds\tilde e_i = \frac{1}{\sqrt{1+\lambda_i^2}}(e_i^h - \lambda_i f_i^v),
        &\ds i=0,1,\dots , n , \\[3ex]
    \ds \tilde n_{\sigma |} =\frac{1}{\sqrt{1+\lambda_\sigma^2}}\big(\lambda_\sigma
    e_\sigma^h +f_\sigma^v \ \big),
    &\ds \sigma=1,\dots , n
    \end{array}
    \end{equation}
    form the orthonormal frames in the tangent space of $\xi(M)$ and in the normal space
    of $\xi(M)$ respectively.
\end{lemma}

    Set
    $$
    (\nabla_X A_\xi)\,Y=\nabla_X(A_\xi Y)-A_\xi \nabla_XY.
    $$
    If we introduce a {\it half tensor } of Riemannian curvature as
    \begin{equation}\label{half}
    r(X,Y)\xi=\nabla_X\nabla_Y\xi-\nabla_{\nabla_XY}\xi,
    \end{equation}
    we can easily  see that
    \begin{equation}\label{A-r}
    -(\nabla_X A_\xi)\,Y=r(X,Y)\xi.
    \end{equation}

    Now, we are able to formulate a basic lemma for our considerations.

\begin{lemma}\cite{Ym1}\label{Form}
    The components of second fundamental form of
    $\xi(M)\subset T_1M$ with respect to the frame (\ref{tang}) are given by
    $$
\begin{array}{ll}
    \tilde \Omega_{\sigma | i j}= &\frac{1}{2}\Lambda_{\sigma i j}
    \Big\{\big< r(e_i,e_j)\xi+  r(e_j,e_i)\xi, f_\sigma \big>+\\[2ex]
    &\hspace{3cm}\lambda_\sigma\left[ \lambda_j \big< R(e_\sigma, e_i)
    \xi, f_j \big> +  \lambda_i  \big<R(e_\sigma, e_j) \xi, f_i \big>\right] \Big\},
\end{array}
   $$
    where  $\Lambda_{\sigma i j}=[(1+\lambda_\sigma^2)(1+\lambda_i^2)(1+\lambda_j^2)]^{-1/2}$
    \ $(i,j=0,1,\dots, n;\,\sigma=1,\dots,n)$.
\end{lemma}

\begin{remark}
    We say that the given unit vector field is {\it holonomic} if $\xi$ is a field of unit
    normals for a family of hypersurfaces in $M$ and {\it non-holonomic} otherwise.
    If the integral trajectories of $\xi$ are geodesics in $M$ then $\xi$ is called
    {\it geodesic} vector field. Evidently, in the case of holonomic geodesic unit
    vector field, $A_\xi$ becomes a usual shape operator for each hypersurface.
    In this case $A_\xi$ is self-adjoint (symmetric), i.e.
    $$
    \big<A_\xi X,Y\big>=\big<X,A_\xi Y\big> \mbox{ and thus } A_\xi^*=A_\xi
    $$
    with respect to some orthonormal frame.  Let
    $R(X,Y)\,\xi=\big[\nabla_X,\nabla_Y\big]\xi-\nabla_{[X,Y]}\,\xi$ be a
    curvature tensor of $M$.
    The non-holonomic shape operator satisfies the {\it non-holonomic
    Codazzi} equation
    \begin{equation}\label{Codazzi}
    R(X,Y)\xi=(\nabla_YA_\xi)X-(\nabla_XA_\xi)Y=r(X,Y)\xi-r(Y,X)\xi.
    \end{equation}
    These facts justify the terminology for the operator $A_\xi$ and the tensor $r$.
\end{remark}

\section{Proof of the main result.}

    In this section we will prove of the following theorem.
\begin{theorem}\label{Hopf}
    Let $\xi$ be a unit geodesic vector field on a unit sphere $S^{n+1}$.
    Then $\xi$ is totally geodesic vector field if and only if $n=2m$
    and $\xi$ is a Hopf vector field\footnote{The nesessary part of the result is true if the vector
    field satisfies a condition of covarian normality, that is $A_\xi A_\xi^*=A_\xi^*A_\xi$.}\label{FN}.
\end{theorem}

\vspace{2ex} {\parindent=0pt \textbf{Sufficient part of the proof.} }

 The sufficiency is a consequence of more general considerations involving the properties
of Killing and Sasakian structure vector fields. Remind that the Hopf vector field is
both Killing and Sasakian structure characteristic vector field.

    We begin with Killing vector fields. Let $M$ be $(n+1)$-dimensional a Riemannian manifold
    admitting a Killing vector field $\xi$. Then $A_\xi$ is skew-symmetric
    \begin{equation}\label{kill}
    \ds
    \big<A_\xi X,Y\big>=-\big<X,A_\xi Y\big> \mbox{ and thus } A_\xi^*=-A_\xi
    \end{equation}
    with respect to some orthonormal frame.

    The field $\xi$ is necessarily {\it geodesic} and evidently $A_\xi\xi=0$.
    For the singular frames we have
    $e_0=f_0=\xi$ and $e_\alpha, f_\alpha \in \xi^\perp$ for all $\alpha=1,\dots,n$.
    Moreover, set $2m=\mathop{\rm rank} A_\xi$. There exists an orthonormal frame
    $$
    e_0=\xi,\ e_1,\dots, e_m,\ e_{m+1},\dots, e_{2m}, \ e_{2m+1},\dots,n$$
    such that
    \begin{equation}\label{Canon1}
    \begin{array}{ll}
    \ds
    A_\xi\,e_\alpha=\lambda_\alpha e_{m+\alpha},\quad A_\xi\,e_{m+\alpha}=
    -\lambda_\alpha e_\alpha
    &\ds\mbox{for $\alpha=1,\dots,m$}\\[2ex]
    \ds A_\xi e_0=0, \quad A_\xi e_\alpha=0 &\ds\mbox{for $\alpha=2m+1,\dots,n$}.
    \end{array}
    \end{equation}
    From the definition of the singular frame (\ref{sing}), we see that one may set
    \begin{equation}\label{Canon2}
    \begin{array}{ll}
    \ds f_\alpha=e_{m+\alpha},\quad f_{m+\alpha}=-e_\alpha &\ds\mbox{ for
    $\alpha=1,\dots,m$}\\[2ex]
    \ds f_0=e_0,\quad f_\alpha=e_\alpha &\ds \mbox{ for $\alpha=2m+1,\dots,n$}.
    \end{array}
    \end{equation}

       A unit vector field $\xi$ is called {\it normal} if $ \big<R(X,Y)Z,\xi\big>=0$
    and {\it strongly normal} if $\big<(\nabla_XA)Y,Z\big>=0$ for all $X,Y,Z\in\xi^\perp$
    \cite{GD-V1}. It is evident that each strongly normal vector field
    is always normal. If $\xi$ is a unit {\it Killing} vector field the converse is also
    true.

     Now we can essentially simplify the result of Lemma \ref{Form}.

    \begin{lemma}\label{Killing}
    Let $\xi$ be a Killing vector field on a Riemannian manifold $M^{n+1}$.
    Suppose that $\mathop{\rm rank} A_\xi=2m$.
    Denote $e_\alpha$ ($\alpha=1,\dots,n$) a
    singular frame for the field $\xi$ satisfying (\ref{Canon1}) and (\ref{Canon2}).
    If $\xi$ is (strongly) normal then the non-zero components of a  second fundamental
    form of $\xi(M)\subset T_1M$ are given by
    $$
    \tilde\Omega_{\sigma|\,m+\sigma\, 0}=-
    \tilde\Omega_{m+\sigma|\,\sigma 0}=\frac{1}{2}\frac{K_\sigma(1-K_\sigma)}{1+K_\sigma}
   \quad \sigma=1,\dots,m,
    $$
    where $K_\alpha$ are the
    sectional curvatures of $M^{n+1}$ along the planes $\xi\wedge e_\alpha $.

    \end{lemma}

   \begin{proof}
   Since $\xi$ is a Killing vector field then \cite{GM}
    \begin{equation}\label{Astar}
    \ds
    A_\xi^*A_\xi X=R(X,\xi)\xi
    \end{equation}
    and hence, the right singular vectors for Killing vector field are the Jacobi
    fields along $\xi$- geodesics. Let $e_1,\dots,e_n$ be an orthonormal frame of Jacobi
    fields. Then
    $$
    R(e_\alpha,\xi)\xi=K_\alpha e_\alpha \ \ \alpha=1,\dots,n,
    $$
    where $K_\alpha$ are the sectional curvatures of $M$ along $\xi\wedge e_\alpha$ planes.

        On the other hand, by definition of singular vectors
    $$
    A_\xi^*A_\xi e_\alpha=\lambda_\alpha^2 e_\alpha.
    $$
     Thus $K_\alpha=\lambda_\alpha^2\geq0$.

    Any Killing vector field $\xi$ satisfies
    \begin{equation}\label{Akill}
    \ds
    r(X,Y)\xi\stackrel{def}{=} -(\nabla_XA_\xi) Y=R(X,\xi)Y.
    \end{equation}
    Therefore, $\ds\tilde\Omega_{\sigma|00}=\big<r(e_0,e_0)\xi,f_\sigma\big>=0.$
    Next, from the  relation
    $$
    \lambda_\alpha\big<e_\sigma,f_\alpha\big>=
    \big<e_\sigma,A_\xi e_\alpha\big>=-\big<A_\xi e_\sigma,e_\alpha\big>=
    -\lambda_\sigma\big<f_\sigma,e_\alpha\big>
    $$
    we find $\ds
    \lambda_\alpha\lambda_\sigma\big<e_\sigma,f_\alpha\big>=
    -\lambda_\sigma^2\,\big<f_\sigma,e_\alpha\big>
    $
    and therefore
    $$
    \begin{array}{ll}
    \tilde \Omega_{\sigma | \alpha 0}
    =& \frac{1}{2}\Lambda_{\sigma\alpha 0}
    \Big\{ \big<
    r(e_\alpha, e_0) \xi+ r(e_0, e_\alpha) \xi, f_\sigma \big>
    + \lambda_\sigma\lambda_\alpha \big< R(e_\sigma, e_0)
    \xi, f_\alpha \big> \Big\}=\\[2ex]
    &\frac{1}{2}\Lambda_{\sigma\alpha 0}
    \Big\{\big< R(e_\alpha,\xi)\xi, f_\sigma \big>+\lambda_\sigma\lambda_\alpha
    \big< R(e_\sigma, \xi)\xi, f_\alpha \big> \Big\}=\\[2ex]
    &\frac{1}{2}\Lambda_{\sigma\alpha 0}
    \Big\{\lambda_\alpha^2\big< e_\alpha, f_\sigma \big>+
    \lambda_\sigma^3\lambda_\alpha\big< e_\sigma, f_\alpha \big> \Big\}=\\[2ex]
    &\frac{1}{2}\Lambda_{\sigma\alpha 0}
    \Big\{\lambda_\alpha^2\big< e_\alpha, f_\sigma \big>-
    \lambda_\sigma^4\big< e_\alpha, f_\sigma \big> \Big\}=\\[2ex]
    &\frac{1}{2}\frac{\lambda_\alpha^2-\lambda_\sigma^4}{\sqrt{(1+\lambda_\sigma^2)(1+\lambda_\alpha^2)}}
    \big< e_\alpha, f_\sigma \big>.
    \end{array}
    $$
    Taking into account (\ref{Canon2}), we get
    \begin{equation}\label{Om}
    \ds
    \tilde\Omega_{\sigma|\,m+\sigma\, 0}=-
    \tilde\Omega_{m+\sigma|\,\sigma 0}=\frac{1}{2}\frac{K_\sigma(1-K_\sigma)}{1+K_\sigma}
   \quad \sigma=1,\dots,m.
    \end{equation}

    Since $\xi $ is strongly normal,
    $$
    r(e_\alpha,e_\beta)\xi\sim\xi,\quad R(e_\alpha,e_\beta)\xi=0
    $$
    for all $\alpha,\beta=1,\dots,n$. Therefore for the remain components we have
    $
    \tilde\Omega_{\sigma|\alpha\beta}=0.
    $
    \end{proof}

One of the most important examples of Killing vector fields is a characteristic vector
field of Sasakian  structure. Using Lemma \ref{Killing} we can prove the following.

\begin{lemma}\label{Sasakian}
    Let $M^{2m+1}$ be Sasakian manifold and $\xi$ be a characteristic  vector field. Then
    $\xi(M)$ is totally geodesic in $T_1M$.
\end{lemma}

\begin{proof}
    Let $M^{2m+1}$ be an odd dimensional  manifold admitting a unit vector field $\xi$,
    linear operator $\varphi$ and 1-form $\eta$ such that
        $$
    \varphi^2X=-X+\eta(X)\xi,\quad \varphi\xi=0,\quad \eta(\varphi X)=0,\quad \eta(\xi)=1
    $$
    for any vector field $X$ on $M$.
     A triple $(\varphi, \xi,\eta)$ is called {\it an
    almost contact structure} on $M$ and the manifold is called {\it an  almost contact }
    manifold.

    If the almost contact manifold is Riemannian with metric
    $g(\cdot,\cdot)=\big<\cdot,\cdot\big>$ and
    $$
    \big<\varphi X,\varphi Y\big>=\big<X,Y\big>-\eta(X)\eta(Y),\quad
    \eta(X)=\big<X,\xi\big>
    $$
    for any vector fields $X$ and $Y$ on $M$  then a tetrad $(\varphi, \xi,\eta, g)$ is
    called an almost contact {\it metric} structure and the manifold is called an almost
    contact {\it metric} manifold. The almost contact metric structure is said to be {\it
    contact} if $d\eta (X,Y)=\big<\varphi X,Y\big>$ for all $X,Y$.

    The vector field $\xi$ is called {\it characteristic vector field } of almost
    contact metric structure. This field is always {\it geodesic}. If, in addition,
     $\xi$ is a Killing vector field, then the almost contact metric manifold is
     called K-contact.
    If, moreover, the Riemannian curvature tensor $R$ satisfies
    \begin{equation}\label{sas2}
    R(X,Y)\,\xi=\big<\xi,Y\big>X-\big<\xi,X\big>Y.
    \end{equation}
    for all vector fields $X,Y$ on $M$, then the K-contact manifold is called
    Sasakian.

    In Sasakian manifold
    \begin{equation}\label{sas1}
    \nabla_X\xi=\varphi X,\quad (\nabla_X\varphi
    )Y=\big<\xi,Y\big>X-\big<X,Y\big>\xi=R(X,\xi)Y.
    \end{equation}

    The  property (\ref{sas2}) implies that the sectional curvature of $M$ along the planes
    involving characteristic vector field $\xi$  is equal to $1$ and $\xi$ is a normal unit
    vector field while (\ref{sas1}) means that $\xi$ is strongly normal. So, we may
    apply Lemma \ref{Killing} and immediately get $\Omega_{\sigma|\,\alpha\beta}\equiv0$.
    \end{proof}

    Now the sufficient part of the Theorem \ref{Hopf} follows immediately from Lemma \ref{Sasakian}.
    Remark, that for the  3-dimensional manifold the we can prove a stronger result.

    \begin{theorem}\label{3dim}
    Let $\xi$ be a unit Killing  vector field on 3-dimensional Riemannian manifold $M^3$.
    If $\xi(M^3)$ is totally geodesic in $T_1M^3$ then
     $M^3$ is either Sasakian and $\xi$ is the characteristic vector field or
     $M^3=M^2\times E^1$ metrically and $\xi$ is the unit vector field of Euclidean factor.
    \end{theorem}
    \begin{proof}
    Indeed, (\ref{Om}) is true for any unit Killing vector field.
    If $\xi(M^3)$ is totally geodesic, then there exist a Jacobi frame $e_1,e_2$ in
    $\xi^\perp$ such that the sectional curvatures $K_{\xi\wedge e_1}=K_{\xi\wedge e_2}=K$ and $K$
    satisfies $K(1-K)=0$.  Since $e_1$ and $e_2$ are Jacobi fields along
    $\xi$-geodesics, $\big<R(e_1,\xi)\xi,e_2\big>=0$ and therefore, for any unit $X$ in
    $\xi^\perp$ it is easy to find that $K_{\xi\wedge X}\equiv K$. Thus, all the sectional
    curvatures along 2-planes involving $\xi$ is equal either 0 or 1. The first case means
    that $\xi$ is a parallel vector field on $M^3$. So, $M^3=M^2\times E^1$ metrically  and $\xi$
    is a unit vector field of Euclidean factor. The second case means that $M^3$ is
    K-contact which means in dimension 3 that $M^3$ is Sasakian and $\xi$ is the
    characteristic vector field.
    \end{proof}
\vspace{2ex}

{\parindent=0pt
\textbf{The necessary part of the proof.}
}
    Suppose now that $\xi$ is geodesic vector field generating a totally geodesic
    submanifold in $T_1S^{n+1}$.

    Since $\xi$ is geodesic vector field,
    $$
    r(e_\alpha,\xi)\xi=\nabla_{e_\alpha}\nabla_{\xi}\xi-\nabla_{\nabla_{e_\alpha}\xi}\xi=-
    A_\xi^2e_\alpha.
    $$
    Since the manifold is of constant curvature 1,
    $$
    \begin{array}{l}
    r(\xi,e_\alpha)\xi=R(\xi,e_\alpha)\xi+r(e_\alpha,\xi)\xi=-e_\alpha-A_\xi^2e_\alpha,\\[1ex]
    \lambda_\alpha\lambda_\sigma\big<R(e_\sigma,\xi)\xi,f_\alpha\big>=
    \lambda_\alpha\lambda_\sigma\big<f_\alpha,e_\sigma\big>=\big<A_\xi
    e_\alpha,A_\xi^*f_\sigma\big>=\big<A_\xi^2e_\alpha,f_\sigma\big>.
    \end{array}
    $$
    So we have
    $$
    \begin{array}{ll}
    \tilde\Omega_{\sigma|\,\alpha 0}=&\frac12\Lambda_{\sigma\alpha0}\left\{
    -2\,\big<A_\xi^2e_\alpha,f_\sigma\big>-\big<e_\alpha,f_\sigma\big>+
    \big<A_\xi^2e_\alpha,f_\sigma\big>\right\}=\\[2ex]
    &
    -\frac12\Lambda_{\sigma\alpha0}\big<A_\xi^2e_\alpha+e_\alpha,f_\sigma\big>.
    \end{array}
    $$
    If $\xi$ is a totally geodesic vector field, then
    $$
    \big<A_\xi^2e_\alpha+e_\alpha,f_\sigma\big>=0
    $$
    for all $\alpha,\sigma=1,\dots,n$. This means that
    $$
    A_\xi^2e_\alpha+e_\alpha=0
    $$
    for all $\alpha$ and therefore  each $e_\alpha$ is the eigenvector for the
    operator $A_\xi^2$ corresponding to a common eigenvalue $-1$. Since $A_\xi$ is a
    {\it real} operator, its eigenvalues are $\pm i$. In this case, with respect to some
    orthonormal frame, the matrix of $A_\xi|_{\xi^\perp}$ takes a box-diagonal form with
    boxes of type\footnote{With additional condition of covariant normality (see the footnote on page \pageref{FN})}
    $$
    \left(\begin{array}{cc}
        0&\pm1\\
        \mp1&0
        \end{array}
    \right).
    $$
    This means that $A_\xi+A_\xi^*=0$ and therefore, $\xi$ is a Killing vector field
    on a sphere. Evidently, the dimension has to be  be odd and the field has to be
    Hopf's one.

\textbf{The proof is complete.}

\begin{remark}
    The Hopf vector field is totally geodesic on the {\it unit} sphere
    only. If $S^{2m+1}(r)$ is a sphere of radius $r$, then $\xi$ is not Sasakian
    structure vector field but still unit Killing normal vector field
    satisfying
    $K_\sigma=\frac{1}{r}$ and therefore
    $$
    \begin{array}{l}
     \ds\tilde\Omega_{\sigma|00}\equiv0,\quad
     \tilde\Omega_{\sigma|\alpha\beta}\equiv0,\\[1ex]
     \ds\tilde \Omega_{\sigma|m+\sigma 0}=-\tilde\Omega_{m+\sigma|\sigma 0}=
     \frac{1}{2\,(1+\frac{1}{r^2})^{3/2}}\,\frac{1}{r^2}(1-\frac{1}{r^2})
     \ \quad  (\alpha,\beta,\sigma=1,\dots,m).
     \end{array}
    $$
    Thus, if $r\ne1$ then $\xi$ does not generate a totally geodesic submanifold.
    Nevertheless, it stays minimal.
\end{remark}

\section{Stability of the Hopf vector field. }
    \subsection{General formula and preparations.}
    Let $M^n$ be a submanifold in a Riemannian space $R^m$. Denote by $\nabla^\perp$ a
    covariant derivative in normal bundle connection. Let $e_i$ ($i=1,\dots,n$) be an
    orthonormal frame on $M^n$ and $\eta $ be a normal variation field along $M^n$.
    Denote by $k_i(\eta)$ an $i$-th principal normal curvature of $M^n$ with respect
    to $\eta$ and $K(e_i,\,\eta)$ the sectional curvature of $R^m$ along a plane
    $e_i\wedge \eta$. Then {\it the volume second variation of $M^n$ with respect to
    $\eta$ } is given by \cite{Dusch}
    \begin{equation}\label{dusch}
   \begin{array}{c}
   \ds \delta^2Vol(\eta)=\\[2ex]
   \ds \int\limits_{M^n}\Bigg\{ \sum\limits_{i=1}^n \big<\nabla_{e_i}^\perp
    \eta, \nabla_{e_i}^\perp \eta\big>
    -|\eta|^2\Big(-\sum\limits_{i,j;\,i\ne j }k_i(\eta)k_j(\eta)+
    \sum\limits_{i=1}^n K(e_i,\eta)\Big)\Bigg\}dV.
   \end{array}
   \end{equation}

    To apply the formula (\ref{dusch}) to our case we should find the normal bundle
    connection for the submanifold $\xi(S^{2m+1})\subset T_1S^{2m+1}$.

    \begin{lemma}\label{Nder}
    Let $X$ and $Y$ be arbitrary vector fields on $S^{n+1}$ ($n=2m$). Denote
    $Y^\perp=Y-\big<\xi,Y\big>\xi$. Then, with respect to
    the Hopf unit vector field $\xi$,
    $$
    \bar\nabla_{\ds X_\xi^{\tau}}Y_\xi^{\nu}=\big(\nabla_X Y^\perp\big)_\xi^{\nu}-
    \frac12\big<\xi,X\big>\big(A_\xi Y^\perp\big)_\xi^{\nu},
    $$
    where $\bar\nabla$ means a covariant derivative with respect to the Levi-Civita
    connection on $T_1S^n$.
    \end{lemma}

    \begin{proof}
    With respect to the Levi-Civita connection $\bar\nabla$ on $T_1M$ we have
    \cite{Bx-V3}
    $$
    \begin{array}{ll}
    \ds\bar\nabla_{\ds X^h}Y^h=(\nabla_XY)^h-\frac12\big(R(X,Y)\xi\big)^{t},
    &\ds\bar\nabla_{\ds X^h}Y^{t}=(\nabla_XY)^{t}+\frac12\big(R(\xi,Y)X\big)^h,\\[2ex]
    \ds\bar\nabla_{\ds X^{t}}Y^h=\frac12\big(R(\xi,X)Y\big)^h,
    &\ds\bar\nabla_{\ds X^{t}}Y^{t}=-\big<Y,\xi\big>X^{t}.
    \end{array}
    $$
    Using this formulas, we get
    $$
    \begin{array}{ll}
    \bar\nabla_{\ds X_\xi^\tau}Y_\xi^\nu&=\bar\nabla_{\ds X^h}(A^*_\xi Y)^h +
    \bar\nabla_{\ds X^h}Y^t-\bar\nabla_{\ds(A^*_\xi X)^t}(A^*_\xi Y)^h-
    \bar\nabla_{\ds(A^*_\xi X)^t}Y^t=\\
    &    \begin{array}{l}
    \ds
    \left[\nabla_X(A_\xi Y)+\frac12R(\xi,Y)X-\frac12R(\xi,A_\xi X)A^*_\xi
    Y\right]^h+\\[1ex]
    \ds\left[\nabla_XY-\frac12R(X,A^*_\xi Y)\xi+\big<Y,\xi\big>A_\xi X\right]^t.
    \end{array}
    \end{array}
    $$

    Since $\xi$ is a Killing vector field, (\ref{kill}) and (\ref{Akill}) are fulfilled and
    we have
    $$
    \begin{array}{ll}
    \bar\nabla_{\ds X_\xi^\tau}Y_\xi^\nu &=
    \left[-(\nabla_X A_\xi) Y-A_\xi(\nabla_XY)+\frac12R(\xi,Y)X+\frac12R(\xi,A_\xi X)A_\xi
    Y\right]^h+\\[1ex]
    &\ \, \begin{array}{l}
    \left[\nabla_XY+\frac12R(X,A_\xi Y)\xi+\big<Y,\xi\big>A_\xi X\right]^t=
    \end{array}\\[1ex]
    &\ \, \begin{array}{l}
        \left[R(X,\xi)Y+\frac12R(\xi,Y)X+\frac12R(\xi,A_\xi X)A_\xi
    Y\right]^h+\\[1ex]
    \left[\frac12R(X,A_\xi Y)\xi+\big<Y,\xi\big>A_\xi X\right]^t+
    (\nabla_XY)_\xi^\nu.
    \end{array}
    \end{array}
    $$
    Since $R$ is a curvature tensor of a unit sphere and keeping in mind (\ref{kill}) and
    (\ref{Astar}), we continue
    $$
    \begin{array}{ll}
    \bar\nabla_{\ds X_\xi^\tau}Y_\xi^\nu &=
    \left[\big<\xi,Y\big>X-\frac12\big<X,Y\big>\xi-\frac12\big<\xi,X\big>Y
    +\frac12\big<A_\xi X,A_\xi Y\big>\xi \right]^h+\\[1ex]
    &\left[-\frac12\big<\xi,X\big>A_\xi Y+\big<\xi,Y\big>A_\xi X\right]^t+
    (\nabla_XY)_\xi^\nu= \\[2ex]
    &
    \left[\big<\xi,Y\big>X-\frac12\big<X,Y\big>\xi-\frac12\big<\xi,X\big>Y
    +\frac12\big<X,Y\big>\xi- \right.\\[1ex]
    &\left.\frac12\big<\xi,X\big>\big<\xi,Y\big>\xi\right]^h+
    \left[-\frac12\big<\xi,X\big>A_\xi Y+\big<\xi,Y\big>A_\xi X\right]^t+
    (\nabla_XY)_\xi^\nu= \\[2ex]
    &
    \left[\big<\xi,Y\big>\big(X-\big<\xi,X\big>\xi\big)-\frac12\big<\xi,X\big>\big(Y-\big<\xi,Y\big>\xi\big)
    \right]^h+\\[1ex]
    &\left[-\frac12\big<\xi,X\big>A_\xi Y+\big<\xi,Y\big>A_\xi X\right]^t+
    (\nabla_XY)_\xi^\nu= \\[2ex]
    &
    \left[\big<\xi,Y\big>\big(A^*_\xi A_\xi X)-\frac12\big<\xi,X\big>\big(A^*_\xi A_\xi Y\big)
    \right]^h+\\[1ex]
    &\left[-\frac12\big<\xi,X\big>A_\xi Y+\big<\xi,Y\big>A_\xi X\right]^t+
    (\nabla_XY)_\xi^\nu= \\[3ex]
    & (\nabla_XY)_\xi^\nu+\big<\xi,Y\big>\big(A_\xi X\big)_\xi^\nu -
    \frac12\big<\xi,X\big>\big(A_\xi Y\big)_\xi^\nu.
    \end{array}
    $$

    Consider, now, $\big(\nabla_XY^\perp)^\nu_\xi$.
    $$
    \begin{array}{ll}
    \big(\nabla_XY^\perp)^\nu_\xi=&\big[A^*_\xi\nabla_XY^\perp\big]^h+
    \big[\nabla_XY^\perp\big]^t=\\[1ex]
    &\left[A^*_\xi\nabla_X\big(Y-\big<\xi,Y\big>\xi\big)\right]^h+
    \left[\nabla_X\big(Y-\big<\xi,Y\big>\xi\big)\right]^t=\\[1ex]
    &\left[A^*_\xi\nabla_XY+\big<\xi,Y\big>A^*_\xi A_\xi X\right]^h+
    \left[\nabla_XY+\big<\xi,Y\big>A_\xi X\right]^t=\\[1ex]
    &(\nabla_XY)_\xi^\nu+\big<\xi,Y\big>\big(A_\xi X\big)_\xi^\nu.
    \end{array}
    $$
    Since $A_\xi\,\xi=0$, we see that $A_\xi Y=A_\xi Y^\perp$ .
    Comparing the results, we get
    $$
    \bar\nabla_{\ds X_\xi^\tau}Y_\xi^\nu =\big(\nabla_XY^\perp)^\nu_\xi-
    \frac12\big<\xi,X\big>\big(A_\xi Y^\perp\big)_\xi^\nu
    $$
    what was claimed.
    \end{proof}

\begin{remark}
    As we see, the tangent component of the latter derivative is zero, which gives
    {\it another proof} of totally geodesic property of the Hopf vector field.
\end{remark}

    Since $Y_\xi^\nu=(Y^\perp)^\nu_\xi$, we may consider only the vectors from
    $\xi^\perp$ in all $\xi$-normal lifting operations. So, keep this in mind in what follows.

\begin{corollary}\label{Nconn}
     Let $\bar\nabla^\perp$ denote a covariant derivative in normal
    bundle of $\xi(M)$. If $\xi$ is a Hopf vector field on the unit sphere $S^{n+1}$
    $(n=2m)$ then
    \begin{equation}\label{nconn}
    \ds
    \bar\nabla^\perp_{\ds X^\tau_\xi}Y_\xi^\nu=\big(\nabla_XY)^\nu_\xi-
    \frac12\big<\xi,X\big>\big(A_\xi Y\big)_\xi^\nu=-\big<\xi,X\big>Y^h+2\big(\nabla_XY\big)^t,
    \end{equation}
    for any $X$ and  $Y\in \xi^\perp$.
\end{corollary}

\begin{proof}
    Indeed, if $\xi$ is a Hopf vector field for any $Y,Z$ we have
     $$
     \begin{array}{ll}
     \big<\big<Y_\xi^\nu,Z_\xi^\nu\big>\big>=&\big<A^*_\xi Y,A^*_\xi Z\big>+
     \big<Y,Z\big>-\big<\xi,Z\big>\big<\xi, Y\big>=\\[1ex]
     &2\big<Y,Z\big>-2\big<\xi,Z\big>\big<\xi, Y\big>.
     \end{array}
     $$

     Let $Y,Z\in\xi^\perp$. Then by Lemma \ref{Nder} and (\ref{kill})
     $$
     \begin{array}{l}
     \big<\big<\bar\nabla_{\ds X_\xi^\tau}Y_\xi^\nu,Z_\xi^\nu\big>\big>=
    \big<\big<\big(\nabla_X Y)^\nu_\xi-
    \frac12\big<\xi,X\big>\big(A_\xi Y\big)_\xi^\nu,Z_\xi^\nu\big>\big> =\\[2ex]
    2\big<\nabla_X Y-\frac12\big<\xi,X\big>A_\xi Y,Z\big>-
    2\big<\nabla_X Y,\xi\big>\big<\xi,Z\big>=\\[2ex]
    2 \big< \nabla_X Y-\big<\nabla_X Y,\xi\big>\xi,Z\big>+\big<\xi,X\big>\big<Y,A_\xi
    Z\big>=\\[2ex]
    \big<\big<2\,(\nabla_X Y)^t,Z^t\big>\big>-\big<\big< \,\big<\xi,X\big>Y^h, -(A_\xi Z)^h
    \big>\big>\\[2ex]
    \big<\big<-\big<\xi,X\big>Y^h+2(\nabla_X Y)^t,(A^*_\xi Z)^h +Z^t\big>\big>=\\[2ex]
    \big<\big<-\big<\xi,X\big>Y^h+2(\nabla_X Y)^t,Z_\xi^\nu\big>\big>.
     \end{array}
    $$
    Since $Z$ is arbitrary, we get the result.
    \end{proof}

\subsection{Second variation of the volume for the Hopf vector field.}

     Let  $\xi$ be a given unit vector field on $M^{n+1}$. Then
     the vector field $\tilde\eta=(\eta)_\xi^\nu$ can be considered as a field of
     volume variation for the submanifold $\xi(M^{n+1})$, where
     $\eta$ is an arbitrary vector field in $\xi^\perp$. For the case of the Hopf
     vector field on $S^{n+1}$ $(n=2m)$, we may choose
     \begin{equation}\label{frame}
     \tilde e_i=\frac{1}{||(e_i)_\xi^\tau||}(e_i)_\xi^\tau=\left\{
     \begin{array}{l}
     e_0^h \mbox{ \ for $i=0$ },\\[2ex]
     \frac{1}{\sqrt{2}}\big(e_\alpha^h+f_\alpha^t\big) \mbox{ \ for $\alpha=1,\dots, n$},
     \end{array}\right.
     \end{equation}
      where $e_0,e^1,\dots,e_n$, $f_0, f_1,\dots, f_n$ form the singular bases for the
     $A_\xi$-  operator, as an orthonormal tangent frame of $\xi(S^{n+1})$. Since
     $\xi(S^{n+1})$ is totally geodesic, the Duschek formula (\ref{dusch}) obtains
     the form
     \begin{equation}\label{duscH}
     \delta^2Vol(\tilde\eta)=
    \int\limits_{S^{n+1}}\Bigg\{ \sum\limits_{i=0}^n \|\nabla_{\ds\tilde e_i}^\perp
    \tilde\eta\,\|
    -\|\tilde\eta\,\|^2\sum\limits_{i=0}^n \tilde K(\tilde e_i,\tilde\eta)\Bigg\}dV,
    \end{equation}
     where $\bar\nabla_{\ds\tilde e_i}^\perp$ is given by (\ref{nconn}) and $\tilde K$
    is given by \cite{BY}
    \begin{equation}\label{Sec}
    \begin{array}{rl}
                \tilde K(\tilde X, \tilde Y)&=
     \big< R(X_1,Y_1)Y_1,X_1\big>-\frac{3}{4}\,\big|\,R(X_1,Y_1)\xi\,\big|^2+\\ [2ex]
     & \frac{1}{4}\,\big|\,R(\xi,Y_2)X_1+R(\xi,X_2)Y_1\big|^2+|X_2|^2|Y_2|^2 -\big<X_2,Y_2\big>^2 + \\[2ex]
     & 3\big< R(X_1,Y_1)Y_2,X_2\big>-\big< R(\xi,X_2)X_1,R(\xi,Y_2)Y_1\big>+\\ [2ex]
     & \big< (\nabla_{X_1}R)(\xi,Y_2)Y_1,X_1\big> +\big<
     (\nabla_{Y_1}R)(\xi,X_2)X_1,Y_1\big>,
    \end{array}
    \end{equation}
    for an orthonormal basis $\tilde X=X_1^h+X_2^v$ , $\tilde Y=Y_1^h+Y_2^v$ ( $X_2,Y_2 \in \xi^\perp$)
    of a 2-plane tangent to $T_1S^{n+1}$.

    \begin{lemma}\label{A}
    Let $\xi$ be the Hopf vector field on $S^{n+1}$ ($n=2m$). Let $\eta\in\xi^\perp$
    generates a normal volume variation $\tilde\eta=(\eta)_\xi^\nu$ for the
    submanifold $\xi(S^{n+1})\subset T_1(S^{n+1})$.
        Let $e_0=\xi, e_1,\dots, e_n$ is a right singular frame for the operator $A_\xi.$
    Then, the Duschek formula (\ref{dusch}) can be reduced to the form
    \begin{equation}\label{varForm}
    \delta^2Vol(\tilde\eta)=
    \int\limits_{S^{n+1}}\Bigg\{4|\nabla_{\ds e_0}\eta\,|^2+2\sum_{\alpha=1}^n |\nabla_{\ds
    e_\alpha}\eta\,|^2-\frac{2n-1}{2}\,|\eta\,|^2 \Bigg\}dV,
    \end{equation}

    \end{lemma}
    \begin{proof}
    Denote $\bar\nabla^\perp$ the normal bundle connection for $\xi(S^{n+1})$.
    Prove, first, that
    $$\sum\limits_{i=0}^n \|\bar\nabla_{\ds\tilde e_i}^\perp
    \tilde\eta\,\|=
    4|\nabla_{\ds e_0}\eta\,|^2+2\sum_{\alpha=1}^n |\nabla_{\ds
    e_\alpha}\eta\,|^2-|\eta\,|^2.
    $$
    Applying (\ref{nconn})  and keeping in mind (\ref{frame}), we find
    $$
    \begin{array}{l}
    \ds
    \bar\nabla^\perp_{\ds\tilde e_0}\tilde\eta=-(\eta)^h+2(\nabla_{\ds e_0}\eta)^t\\[1ex]
    \ds
    \bar\nabla^\perp_{\ds\tilde e_\alpha}\tilde\eta=\sqrt{2}\,(\nabla_{\ds e_\alpha}\eta)^t
    \end{array}
    $$
    Therefore
    $$
    \begin{array}{ll}
    \ds
    \sum\limits_{i=0}^n \|\bar\nabla_{\ds\tilde e_i}^\perp \tilde\eta\,\|=&
    |\eta\,|^2+4\,|\nabla_{\ds e_0}\eta\,|^2+2\sum_{\alpha=1}^n\left(\,|\nabla_{\ds
    e_\alpha}\eta\,|^2-\big<\nabla_{\ds e_\alpha}\eta\,,\xi\big>^2\right)=\\[1ex]
    &
    |\eta\,|^2+4\,|\nabla_{\ds e_0}\eta\,|^2+2\sum_{\alpha=1}^n\left(\,|\nabla_{\ds
    e_\alpha}\eta\,|^2-\big<\eta\,,\nabla_{\ds e_\alpha}\xi\big>^2\right)=\\[1ex]
    &
    |\eta\,|^2+4\,|\nabla_{\ds e_0}\eta\,|^2+2\sum_{\alpha=1}^n\,|\nabla_{\ds
    e_\alpha}\eta\,|^2-2\sum_{\alpha=1}^n\big<\eta\,,f_\alpha\big>^2=\\[1ex]
    &
    |\eta\,|^2+4\,|\nabla_{\ds e_0}\eta\,|^2+2\sum_{\alpha=1}^n\,|\nabla_{\ds
    e_\alpha}\eta\,|^2-2\,|\eta\,|^2=\\[1ex]
    &
    4\,|\nabla_{\ds e_0}\eta\,|^2+2\sum_{\alpha=1}^n\,|\nabla_{\ds
    e_\alpha}\eta\,|^2- |\eta\,|^2
    \end{array}
    $$
    and the proof is complete.

    Now prove that
    $$
    \|\tilde \eta\,\|^2\,\sum_{i=0}^n\tilde K(\tilde e_i,\tilde \eta\,)=
    \frac12|\eta\,|^2+(n-2)\,|\eta\,|^2.
    $$
    The sectional curvature of $T_1M^n$ is given by (\ref{Sec}) in presumption that
    $\tilde X$ and $\tilde Y$ are orthonormal. To find
    $\|\tilde \eta\,\|^2 \tilde K(\tilde e_i,\tilde \eta\,)$
    we can use (\ref{Sec}) setting $\tilde Y=\tilde \eta$ and keeping in mind that
    $\tilde \eta$ is of arbitrary length.

    Now set $\zeta=A_\xi\eta$ and then $\tilde \eta=\zeta^h+\eta^t.$
    Evidently,
    \begin{equation}\label{hint}
    \big<\eta,\zeta\big>=0,\quad\big<\eta,\xi\big>=\big<\zeta,\xi\big>=0,
    \quad |\eta\,|=|\zeta|.
    \end{equation}

    Set $\tilde X=\tilde e_0=e_0^h$ and $\tilde Y=\zeta^h+\eta^t$. Then
    $$
    \begin{array}{ll}
    \ds \|\tilde\eta\,\|^2
    \tilde K(\tilde e_0,\tilde \eta\,)=&\ds\big< R(e_0,\zeta)\zeta,e_0\big>-\frac{3}{4}\,\big|\,R(e_0,\zeta)\xi\,\big|^2+
    \frac{1}{4}\,\big|\,R(\xi,\eta)e_0\,\big|^2=\\[1ex]
    &\ds |\zeta|^2-\frac34\,|\zeta|^2+\frac14\,|\eta\,|^2=\frac12\,|\eta\,|^2
    \end{array}
    $$
    To find $\tilde K(\tilde e_\alpha,\tilde \eta\,)$ for $\alpha=1,\dots,n$ we set
    $\tilde X=\frac{1}{\sqrt{2}}(e_\alpha^h+f_\alpha^t)$ and therefore
    $$
    X_1=\frac{1}{\sqrt{2}}e_\alpha, \ X_2=\frac{1}{\sqrt{2}}f_\alpha,\ Y_1=\zeta, \
    Y_2=\eta
    $$
    in application of (\ref{Sec}). Thus we have
    $$
    \begin{array}{rl}
    \ds \|\tilde\eta\,\|^2
    \tilde K(\tilde e_\alpha, \tilde \eta)&=\ds
     \frac{1}{2}\,\big< R(e_\alpha,\zeta)\zeta,e_\alpha\big>-\frac{3}{8}\,\big|\,R(e_\alpha,\zeta)\xi\,\big|^2+\\ [2ex]
     &\ds \frac{1}{8}\,\big|\,R(\xi,\eta)e_\alpha+R(\xi,f_\alpha)\zeta\big|^2+\frac{1}{2}|\eta\,|^2 -\frac{1}{2}\big<f_\alpha,\eta\big>^2 + \\[2ex]
     & \ds \frac32\big< R(e_\alpha,\zeta)\eta,f_\alpha\big>-\frac12\big<
     R(\xi,f_\alpha)e_\alpha,R(\xi,\eta)\zeta\big>.
    \end{array}
    $$
    Now set $\eta^\alpha=\big<\eta,f_\alpha\big>$. Then
    $\zeta^\alpha=\big<\zeta,e_\alpha\big>=-\big<\eta,f_\alpha\big>=-\eta^\alpha.$
     Keeping in mind (\ref{hint}) we get
    $$
    \begin{array}{rl}
     \|\tilde\eta\,\|^2
    \tilde K(\tilde e_\alpha, \tilde \eta)=&\ds
     \frac{1}{2}\,|\eta\,|^2-\frac12(\eta^{\,\alpha})^2+
     \frac{1}{8}\,\big|\big<e_\alpha,\eta\,\big>+
     \big<f_\alpha,\zeta\,\big>\big|^2+
     \frac{1}{2}|\eta\,|^2 - \\[2ex]
     &\frac{1}{2}(\eta^{\,\alpha})^2 -
      \frac32\big<e_\alpha,\eta\,\big>\big<f_\alpha,\zeta\,\big>=
     |\eta\,|^2- (\eta^{\,\alpha})^2-\big<\eta,e_\alpha\big>^2.
    \end{array}
    $$
    It is also easy to see that
    $$
    \sum_{\alpha=1}^n (\eta^{\,\alpha})^2=\sum_{\alpha=1}^n
    \big<\eta,e_\alpha\big>^2=|\eta\,|^2.
    $$
    Hence we have
    $$
    \|\tilde\eta\,\|^2\left(
     \tilde K(\tilde e_0,\tilde\eta)+\sum_{\alpha=0}^n \tilde K(\tilde e_\alpha,\tilde\eta)
     \right)=
     \frac12 \,|\eta\,|^2+(n-2)\,|\eta\,|^2.
    $$
    Combining the results we get what was claimed.
    \end{proof}

    \begin{proposition}\label{Hopf3}
    The Hopf vector field on the unit 3-sphere is stable.
    \end{proposition}

    \begin{proof}
    Since the subspace of right (and left) singular frames for the Hopf
    vector field coincides with $\xi^\perp$, we may choose the other two Hopf vector
    fields on $S^3$ as $e_1$ and $e_2$. Then
    \begin{equation}\label{proiz}
    \begin{array}{ll}
    \nabla_{\ds e_1}e_1=0, & \nabla_{\ds e_1}e_2=-e_0, \\
     \nabla_{\ds e_2}e_1=e_0, & \nabla_{\ds e_2}e_2=0.
    \end{array}
    \end{equation}
    Set $\eta=\eta^{1}e_1+\eta^2e_2$. Set $|{\rm grad}\,
    \eta^{\,\alpha}|^2=\big(e_1(\eta^{\,\alpha})\big)^2+\big(e_2(\eta^{\,\alpha})\big)^2$
    ($\alpha=1,2$). Then, using (\ref{proiz}) we find
    $$
    \sum_{\alpha=1}^2 |\nabla_{\ds e_\alpha}\eta\,|^2=\sum_{\alpha =1}^n |\,{\rm grad}\,
    \eta^{\,\alpha}|^2+|\eta\,|^2.
    $$
    Therefore,
    $$
     \delta^2Vol(\tilde\eta)=
    \int\limits_{S^{n+1}}\Bigg\{4|\nabla_{\ds e_0}\eta\,|^2+2\sum_{\alpha=1}^n |{\rm grad}\,
    \eta^{\,\alpha}|^2+\frac{1}{2}\,|\eta\,|^2 \Bigg\}dV >0
    $$
    which means that $\xi(S^3)$ is stable.
    \end{proof}

    \begin{proposition}\label{Hopfn}
    The Hopf vector field on the unit $n$-sphere for $n=2m+1>3$ is unstable.
    \end{proposition}
    \begin{proof}
    Choose the vectors of the singular frame such that $e_0=\xi$ while the other
    $e_1,\dots,e_{2m}$ are the horizontal lifts of vectors of orthonormal frame $q_k$
    of $CP^m$ with respect to the Hopf fibration
    $S^{2m+1}\stackrel{S^1}{\longrightarrow} CP^m$. Then $\nabla_{ q_k}q_j=0$ and
    $Jq_{2k}=q_{2k+1}$ ($k,j=1,\dots,m$) for complex structure of $CP^n$ (see \cite{HY}). Then along the
    fiber, i.e. along the integral curves of $e_0=\xi$, the following table of
    non-zero covariant derivatives can be achieved \cite{HY,J}.
    $$
    \begin{array}{ll}
    \nabla_{e_0}e_{2k}=e_{2k-1}, & \nabla_{e_0}e_{2k-1}=-e_{2k}, \\[1ex]
    \nabla_{ e_{2k}}e_{0}=e_{2k-1}, & \nabla_{ e_{2k-1}}e_{0}=-e_{2k}\\[1ex]
    \nabla_{ e_{2k-1}}e_{2k}=e_{0}, & \nabla_{ e_{2k}}e_{2k-1}=-e_{0},
    \end{array}
    $$
    where $k=1,\dots,m$.

    Let $\eta=\eta^{\,2k-1}e_{2k-1}+\eta^{\,2k}e_{2k}$ be a variation field. Then,
    using a table of derivatives, we find
    $$
    \begin{array}{l}
    \nabla_{e_0}\,\eta=\big(e_0(\eta^{\,2k-1})+\eta^{\,2k}\big)\,e_{2k-1}+
                    \big(e_0(\eta^{\,2k})-\eta^{\,2k-1}\big)\,e_{2k-1} \\[1ex]
    \nabla_{e_{2t-1}}\,\eta=e_{2t-1}(\eta^{\,2k-1})\,e_{2k-1}+
                    e_{2t-1}(\eta^{\,2k})\,e_{2k}+ \delta_{st}\,\eta^{\,2k}\,e_0 \\[1ex]
    \nabla_{e_{2t}}\,\eta=e_{2t}(\eta^{\,2k-1})\,e_{2k-1}+
                  e_{2t}(\eta^{\,2k})\,e_{2k}- \delta_{kt}\,\eta^{\,2k-1}\,e_0 ,
    \end{array}
    $$
    where $\delta_{kt}$ is the Kronecker symbol.

    Set
    $$
    |{\rm grad}\,\eta^{\,\sigma}|^2=\sum_{\alpha=1}^n \big[e_\alpha(\eta^{\,\sigma})\big]^2.
    $$
    Then
    $$
    \sum_{\alpha=1}^n |\nabla_{e_\alpha}\eta\,|^2=\sum_{\sigma=1}^n |{\rm
    grad}\,\eta^{\,\sigma}|^2+ |\eta\,|^2
    $$
    and
    $$
    |\nabla_{e_0}\eta\,|^2=\sum_{k=1}^m \left[
    \big(e_0(\eta^{\,2k})-\eta^{\,2k-1}\big)^2+\big(e_0(\eta^{\,2k-1})+\eta^{\,2k}\big)^2\right].
    $$
    To prove the instability, we should find variation field $\eta$ providing a
    negative sign for the second volume variation. So, choose
    $$
    \eta= \cos t e_{2k-1}+\sin t e_{2k},
    $$
    where $t$ is an arc-length parameter on $e_0$-curves. Then
    $$
    \nabla_{e_0}\eta\,=0,\quad {\rm grad}\,\eta^{\,\sigma}=0
    $$
    and for the integrand in the Duscheck formula (\ref{varForm}) we have
    $$
    \frac{5-2n}{2}|\eta\,|^2<0
    $$
    for $n>2$, which completes the proof.

    \end{proof}

\section{ Sectional curvature of the Hopf vector field.}

Since the Hopf vector field is totally geodesic in $T_1S^{n+1}$, the sectional
curvature of $\xi(S^{n+1})$ is completely defined by the curvature of $T_1S^{n+1}$
along the planes, tangent to $\xi(S^{n+1})$. We can easily find it applying
(\ref{lifts}) and (\ref{Sec}) to the Hopf vector field. Remind that the Hopf vector
field is a characteristic one for the contact metric structure on the spheres. In
contact metric geometry, the sections, containing characteristic vector field, are
called $\xi$-sections, while the sections of type $X\wedge \varphi X$ for
$X\in\xi^\perp$ are called $\varphi$-sections. Using the notion of $\xi$-tangential
lift (\ref{lifts}), for any $X\in \xi^\perp$ we call a 2-plane $X^\tau_\xi\wedge
(\xi)_\xi^\tau$ as a $\xi$-tangential lift of $\xi$-section and a 2-plane
$X^\tau_\xi\wedge(\nabla_X\xi)^\tau_\xi,$ as a $\xi$-tangential lift of
$\varphi$-section. The following assertion holds.

\begin{theorem}\label{Cor}
    The sectional curvature of $\xi(S^{2m+1})$ for the Hopf vector field varies between
$\frac14$ and $\frac54.$ The curvature is minimal for $\xi$-tangential lift of
$\xi$-section and maximal for $\xi$-tangential lift of $\varphi$-section
    \end{theorem}
The proof is elementary consequence of the Proposition below. It should be mentioned
that the curvature of $T_1S^{n+1}$ varies between $0$ and $5/4$ \cite{Ym2}. The
theorem clarifies geometrical meaning of the maximal curvature. The minimal curvature
is, geometrically, the $\xi$-sectional curvature of natural Sasakian structure  on
$T_1S^{n+1}$ (which,as well known, equals to $1/4$ after rescaling).

\begin{proposition} Let $\xi$ be the Hopf vector field on $S^{n+1}$ ($n=2m$).
Let $X_\xi^\tau$, $Y_\xi^\tau$ be $\xi$-tangential lifts of orthonormal vectors $X$
and $Y$ respectively. The sectional curvature $\tilde K(X_\xi^\tau, Y_\xi^\tau)$ of
$\xi(S^{n+1})$ along the 2-plane ($X_\xi^\tau$, $Y_\xi^\tau$) is given by
    $$
    \tilde K(X_\xi^\tau,Y_\xi^\tau)=
    \frac{\ds 1-\frac34\big[\,\big<\xi,X\big>^2+\big<\xi,Y\big>^2\,\big]+
            \frac32\,\big<A_\xi X,Y\big>^2}{\ds 2-\big[\,\big<\xi,X\big>^2+\big<\xi,Y\big>^2\,\big]}
    $$
\end{proposition}

    \begin{proof}
    Let $X$ and $Y$ be {\it unit mutually orthogonal} vector fields on $S^{n+1}$. Then
    $X_\xi^\tau=X^h-(A_\xi X)^t$ and $Y_\xi^\tau=Y^h-(A_\xi Y)^t$ form a basis of
    elementary 2-plane, tangent to $\xi(S^{n+1})$. This basis is not
    orthonormal, since
    $$
    \begin{array}{l}
    \ds
    \|X_\xi^\tau\|^2=|X|^2+|A_\xi X|^2=2-\big<\xi,X\big>^2,\\[1ex]
    \ds \|Y_\xi^\tau\|^2=|Y|^2+|A_\xi Y|^2=2-\big<\xi,Y\big>^2, \\[1.5ex]
    \ds \big<\big<X_\xi^\tau,Y_\xi^\tau\big>\big>=\big<X,Y\big>+\big<A_\xi X,A_\xi
    Y\big>=-\big<\xi,X\big>\big<\xi,Y\big>.
    \end{array}
    $$
    Therefore, the norm of bivector $X_\xi^\tau\wedge Y_\xi^\tau$ is
    \begin{equation}\label{bivec}
    \|X_\xi^\tau\wedge Y_\xi^\tau\|=4-2\,\big<\xi,X\big>^2-2\,\big<\xi,Y\big>^2.
    \end{equation}
    Now set $X_1=X,\ Y_1=Y,\ X_2=A_\xi X,\ Y_2= A_\xi Y$ and apply (\ref{Sec}). We
    have
    $$
    \begin{array}{l}
    \big< R(X_1,Y_1)Y_1,X_1\big>=1, \\[2ex]
    \big|\,R(X_1,Y_1)\xi\,\big|^2=\big|\big<\xi,Y\big>\,X-\big<\xi,X\big>\,Y\big|^2=
    \big<\xi,Y\big>^2+\big<\xi,Y\big>^2, \\[2ex]
    \big|\,R(\xi,Y_2)X_1+R(\xi,X_2)Y_1\big|^2=\big|\big<A_\xi Y,X\big>\xi-\big<\xi,X\big>A_\xi Y+\\[1ex]
    \ \big<A_\xi X,Y\big>\xi-\big<\xi,Y\big>A_\xi X\big|^2=
    \big<\xi,X\big>^2|A_\xi Y|^2+\big<\xi,Y\big>^2|A_\xi X|^2+ \\[1ex]
    \ 2\,\big<\xi,X\big>\big<\xi,Y\big>\big<A_\xi X,A_\xi Y\big>=
    \big<\xi,X\big>^2(1-\big<\xi,Y\big>^2\,)+\big<\xi,Y\big>^2(1-\big<\xi,X\big>^2\,)-\\[1ex]
    \  2\,\big<\xi,X\big>^2\big<\xi,Y\big>^2=\big<\xi,X\big>^2+\big<\xi,Y\big>^2-
    4\,\big<\xi,X\big>^2\big<\xi,Y\big>^2,
     \\[2ex]
    |X_2|^2|Y_2|^2 -\big<X_2,Y_2\big>^2=|A_\xi X|^2|A_\xi Y|^2- \big<A_\xi X,A_\xi
    Y\big>^2=\\[1ex]
    \ (1-\big<\xi,X\big>^2\,)(1-\big<\xi,Y\big>^2\,)-\big<\xi,X\big>^2\big<\xi,Y\big>^2=
    1-\big<\xi,X\big>^2-\big<\xi,Y\big>^2,
    \\[2ex]
    \big< R(X_1,Y_1)Y_2,X_2\big>=\big<-\big<A_\xi Y,X\big>\,Y,A_\xi X\big>=\big<A_\xi
    X,Y\big>^2, \\[2ex]
    \big< R(\xi,X_2)X_1,R(\xi,Y_2)Y_1\big>=\big<\,\big<\xi,X\big>\,A_\xi X,\big<\xi,Y\,\big>\,A_\xi
    Y\big>=-\big<\xi,X\big>^2\big<\xi,Y\big>^2.
    \end{array}
    $$
    Substituting the latter equalities into (\ref{Sec}) and dividing the result by
    (\ref{bivec}), we get
    $$
    \tilde K(X_\xi^\tau,Y_\xi^\tau)=
    \frac{\ds 1-\frac34\big[\,\big<\xi,X\big>^2+\big<\xi,Y\big>^2\,\big]+
            \frac32\,\big<A_\xi X,Y\big>^2}{\ds 2-\big[\,\big<\xi,X\big>^2+\big<\xi,Y\big>^2\,\big]}
   $$
    \end{proof}

\vspace{1cm}

\noindent
Department of Geometry,\\
Faculty of Mechanics and Mathematics,\\
Kharkiv National University,\\
Svobody Sq. 4,\\
 61077, Kharkiv,\\
Ukraine.\\
e-mail: yamp@univer.kharkov.ua

\end{document}